\documentclass[12pt]{amsart}

\newtheorem{thm}{Theorem}
\newtheorem{prop}{Proposition}
\newtheorem{defin}{Definition}

\newtheorem{corol}{Corollary}

\newcommand{\N}{{\mathbb N}}

\newcommand{\R}{{\mathbb R}}\newcommand{\RR}{{\mathbb R}^2}

\newcommand{\Z}{{\mathbb Z}}

\newcommand{\bo}{\partial} 

\newcommand{\inter}{\text{interior}}     

\newcommand{\tz}{\tilde{z}}

\newcommand{\al}{\alpha}

\newcommand{\ga}{\gamma}\newcommand{\Ga}{\Gamma}\newcommand{\tga}{\tilde{\gamma}}
\newcommand{\de}{\delta}\newcommand{\tde}{\tilde{\delta}}

\newcommand{\vap}{\varphi}
\newcommand{\si}{\sigma}
\newcommand{\het}{\theta}
\newcommand{\Om}{\Omega}
\newcommand{\ups}{\upsilon}

\begin{document}

\bibliographystyle{plain}

\title[billiard and rotation numbers]
{The billiard ball problem and rotation numbers}

\author{Eugene Gutkin}

\address{UMK and IMPAN \\
Chopina 12/18\\
87 -- 100 Torun\\
Poland}

\email{gutkin@mat.uni.torun.pl,\ gutkin@impan.pl}

\keywords{billiard map, liouville measure, canonical involution, rotation vector, birkhoff ergodic theorem, rotation number}

\subjclass{37E45, 37A10}
\date{\today}

\begin{abstract}
We introduce the concepts of rotation numbers and rotation vectors
for billiard maps. Our approach is based on the birkhoff ergodic
theorem. We anticipate that it will be useful, in particular, for
the purpose of establishing the non-ergodicity of billiard in
certain domains.
\end{abstract}

\maketitle

\section{Motivation}       \label{intro}
%
In \cite{Po82} Poincare introduced and studied the concept of
rotation number for orientation preserving circle homeomorphisms.
Since then, rotation numbers and their generalizations became
instrumental in several areas of mathematics. The purpose of a
brief outline that follows is to give the reader an idea about the
spread of the broadly interpreted concept of rotation number in
the contemporary mathematical literature.

\medskip

E. Ghys \cite{Gh01}  used rotation numbers to show that certain
classes of groups cannot act on the circle. A.I. Stern
\cite{Sht07} encountered rotation numbers in the framework of
quasi-representations of groups. M. Misiurewicz and others
extended the framework of rotation numbers in several directions.
In particular, rotation numbers become rotation sets when circle
homeomorphisms are replaced by arbitrary continuous self-mappings
of the circle or by homeomorphisms of multi-dimensional tori. We
recommend the text \cite{Mi08} for a comprehensive exposition and
an extensive bibliography.

\medskip

The concept of rotation number is quite flexible which explains the differences
between approaches to this concept in various dynamics frameworks \cite{Ru85}.
Here we are interested in the uses of rotation numbers for billiard dynamics.
But even in the subject of billiards there are several versions of the notion of rotation numbers.
The work \cite{BlMiSh06} defines rotation numbers in the framework of billiard flows.
In the interpretation of \cite{BlMiSh06}, rotation numbers provide information about
the way orbits of a billiard flow wind around an obstacle inside a billiard table.
This approach is especially useful for studying periodic billiard orbits.
However, it is designed and developed for rather special billiard tables.

\medskip

In the present work we introduce a certain analog of the rotation number for
simply connected billiard tables. If the billiard table is multiply connected, then
we develop the concept of rotation vectors. We use the billiard map, as opposed to
the billiard flow; we rely on the birkhoff ergodic theorem to define the
rotation vectors. Thus, in general, our approach leads to vector valued functions
defined for almost all phase points. Our rotation vectors
make sense for all of the basic paradigms in billiard dynamics: {\em elliptic, hyperbolic and
parabolic} \cite{Gut96,Gut03}. For simplicity of exposition, we restrict out attention
to the planar billiard dynamics. It is straightforward to extend our treatment to
billiard domains on surfaces of any constant curvature \cite{GSG99}.



\medskip

\medskip

\noindent{\bf Acknowledgements:} The work was completed during a visit to UCLA in
January -- February 2009. It is a pleasure to thank the Department of Mathematics
at UCLA for its hospitality.

\vspace{4mm}

\section{The billiard ball problem for simply connected billiard tables}   \label{simply}
Let $\Om\subset\RR$ be a bounded, connected domain with a
piecewise smooth boundary $\bo\Om$. Let $\Phi=\Phi(\Om)$ be the
{\em billiard map phase space}. The {\em billiard map}
$\vap:\Phi\to\Phi$ has a natural invariant measure, the {\em
liouville measure}. Let $(s,\het)$ be the standard coordinates in
$\Phi$ \cite{Gut96}. They satisfy $0\le s \le |\bo\Om|$, where
$|\bo\Om|$ is the perimeter of the billiard table, and $0\le \het
\le \pi$. The liouville measure has a density: $d\nu = \sin\het
dsd\het$. Thus $\nu(\Phi)=2|\bo\Om|$. We refer to \cite{Gut96} and
\cite{Gut03} for details. There is a natural projection
$p:\Phi\to\bo\Om$. Let $z\in\Phi$ be a phase point. Then
$p(z)\in\bo\Om$ is the  {\em footpoint} of $z$.

\vspace{4mm}
From now and until the end of this section we assume that $\Om$ is simply connected.
Thus, $\bo\Om$ is a simple closed curve.
Choosing a reference point in $\bo\Om$, endowing $\bo\Om$ with the positive orientation,
and parameterizing $\bo\Om$ by the arclength, we identify
$\bo\Om$ with the quotient $\R/|\bo\Om|\Z$. We will regard $\R/|\bo\Om|\Z$ as a circle of perimeter
$|\bo\Om|$. Let $x_1,x_2\in \bo\Om$ be arbitrary.
Let $0\le s_1,s_2< |\bo\Om|$ be their arclength coordinates.
Then there is a unique
number $0\le \xi=\xi(x_1,x_2)<|\bo\Om|$ such that either $s_1+\xi=s_2$ or $s_1+\xi=s_2+|\bo\Om|$.

Let $z\in\Phi$. We define the {\em footpoint increment function} $\xi:\Phi\to\R$ via
\begin{equation}   \label{increm_eq}
\xi(z)=\xi(p(z),p(\vap(z))).
\end{equation}
If the limit
\begin{equation}   \label{rot_num_eq}
\upsilon(z)=\lim_{n\to\infty}\frac{1}{n}[\sum_{k=0}^{n-1}\xi(\vap^k(z)]
\end{equation}
exists, we say that $\upsilon(z)$ is the {\em footpoint gain} of the phase point $z\in\Phi$.
The birkhoff ergodic theorem immediately implies the following.

\medskip

\begin{prop}    \label{foot_gain_prop}
The footpoint gain $\upsilon(z)$ exists for almost all $z\in\Phi$ with respect to the
liouville measure. The function $\upsilon:\Phi\to\R$ is measurable,
invariant under $\vap$, and satisfies  $0\le \upsilon(\cdot) \le |\bo\Om|$.
We have
$$
0\le \frac{\int_{\Phi}\upsilon(z)d\nu(z)}{\nu(\Phi)} \le |\bo\Om|.
$$
\end{prop}

\medskip

It is often convenient to normalize the function $\upsilon(z)$ and define
the {\em rotation number function} $\rho(z)=\upsilon(z)/|\bo\Om|$.
We reformulate Proposition~~\ref{foot_gain_prop}
in terms of the rotation number.

\begin{prop}    \label{rot_num_prop}
The rotation number $\rho(z)$ exists for almost all $z\in\Phi$ with respect to the
liouville measure. The function $\rho:\Phi\to\R$ is measurable,
invariant under $\vap$, and satisfies  $0\le \rho(\cdot) \le 1$.
The average rotation number satisfies
\begin{equation}   \label{av_rot_num_eq}
0\le \frac{\int_{\Phi}\rho(z)d\nu(z)}{\nu(\Phi)} \le 1.
\end{equation}
\end{prop}
\begin{thm}    \label{rot_num_thm}
Let $\Om\subset\RR$ be a bounded, connected, simply connected domain with a piecewise smooth boundary.
Let $\vap:\Phi\to\Phi$ be the billiard map; let $\rho:\Phi\to\R$ be the rotation number
function.

\noindent 1. Let $X\subset[0,1]$ be a measurable set. Denote by $1-X\subset[0,1]$ the set
of numbers $y=1-x$, where $x\in X$. Then $\rho^{-1}(X)$ and $\rho^{-1}(1-X)$ are measurable
subsets of $\Phi$, and
$\nu(\rho^{-1}(X))=\nu(\rho^{-1}(1-X))$.

\noindent 2. We have
\begin{equation}    \label{mean_rot_eq}
\frac{\int_{\Phi}\rho(z)d\nu(z)}{\nu(\Phi)} = \frac12.
\end{equation}

\noindent 3. If the billiard in $\Om$ is ergodic then $\rho(z)=1/2$ for almost all $z$.
\begin{proof}
The space $\Phi$ carries two canonical involutions: $\si:\Phi\to\Phi$ and $\tau:\Phi\to\Phi$ \cite{Gut96}.
Let $z=(s,\het)\in\Phi$. Then $\si(z)=(s,\pi-\het)$. In order to define $\tau$, we identify $\Phi$ with
the set of oriented line segments in $\Om$ whose endpoints belong to $\bo\Om$. Then $\tau$ corresponds
to the direction reversal for the segments. Both involutions preserve the liouville measure; we have
$$
\vap=\si\circ\tau.
$$
Therefore
\begin{equation}    \label{invol_eq}
\si\circ\vap\circ\si=\tau\circ\vap\circ\tau=\vap^{-1},
\end{equation}
i. e., both $\si$ and $\tau$ conjugate the billiard map with its inverse. By definition of the footpoint increment,
we have
$$
\xi(\tau(z))=|\bo\Om|-\xi(z).
$$
This observation plus equation~~\eqref{invol_eq} imply the relationship
\begin{equation}    \label{revers_eq}
\upsilon(\tau(z))=|\bo\Om|-\upsilon(z).
\end{equation}
Equation~~\eqref{revers_eq} means, in particular, that the limits in equation~~\eqref{rot_num_eq}
defining $\upsilon(z),\,\upsilon(\tau(z))$ exist or do not exist simultaneously. We rewrite
equation~~\eqref{revers_eq} as a relationship for rotation numbers:
\begin{equation}    \label{revers_eq1}
\rho(\tau(z))=1-\rho(z).
\end{equation}
\medskip

Equation~~\eqref{revers_eq1} immediately implies the first claim. The second claim is
a straightforward consequence of the first. The third claim follows from the second and the
birkhoff ergodic theorem.
\end{proof}
\end{thm}
\medskip

\vspace{5mm}

%

\medskip


\medskip

\medskip

\vspace{3mm}

\section{The billiard ball problem for multiply connected billiard tables}   \label{multi}
In this section we study the billiard dynamics in a {\em multiply connected planar domain}.
We will use the setting and the notation established in the opening paragraph of section~~\ref{simply}.
We will say that the
domain $\Om$ is $q$-connected if its boundary  is a disjoint union
of $q$ connected components:
$$
\bo\Om=\bo_1\Om\cup\cdots\cup\bo_q\Om.
$$
Thus
$$
|\bo\Om|=|\bo_1\Om|+\cdots+|\bo_q\Om|.
$$
We choose a reference point on each component,
$\bo_{\al}\Om,1\le\al\le q$. We parameterize each component by the
arclength with respect to the positive orientation. Thus, we
identify  each curve $\bo_{\al}\Om$ with the circle
$\R/|\bo_{\al}\Om|\Z$.

\medskip

Let $z\in\Phi$. The footpoints $p(z),p(\vap(z))$ may belong to different components of $\bo\Om$.
In view of this observation, there is no counterpart of the footpoint increment function $\xi(z)$
of equation~~\eqref{increm_eq}. We will directly define the counterpart of the footpoint gain
$\upsilon(z)$ in  equation~~\eqref{rot_num_eq}.

\medskip

We fix a boundary component $\bo_{\al}\Om$. Let $z\in\Phi$, and let $N\in\N$. Let
$0\le k_0<\cdots<k_n\le N$ be the consecutive times such that $p(\vap^{k_i}(z))\in\bo_{\al}\Om$.
Let $0\le s_0,\dots,s_n<|\bo_{\al}\Om|$ be their coordinates.
For every $0\le i \le n-1$ there is a unique number
$0\le \xi_i<|\bo_{\al}\Om|$ such that either
$s_i+\xi_i=s_{i+1}$ or $s_i+\xi_i=s_{i+1}+|\bo_{\al}\Om|$.
Suppose that the limit
\begin{equation}   \label{comp_rot_num_eq}
\upsilon_{\al}(z)=\lim_{N\to\infty}\frac{\xi_0+\cdots+\xi_{n-1}}{N}
\end{equation}
exists. Then $\upsilon_{\al}(z)$ is the $\al$-component of the {\em footpoint gain  vector}
%
\begin{equation}   \label{gain_vec_eq}
\vec{\upsilon}(z)=(\upsilon_1(z),\dots,\upsilon_q(z)).
\end{equation}

\medskip

Let $\vec{a}=(a_1,\dots,a_q),\,\vec{b}=(b_1,\dots,b_q)\in\R^q$ be arbitrary vectors.
We will use the notation $\vec{a}\prec\vec{b}$ to mean that $a_{\al}\le b_{\al}$ for $1\le \al \le q$.

The birkhoff ergodic theorem immediately implies the following.

\begin{prop}    \label{vec_foot_gain_prop}
The footpoint gain vector $\vec{\upsilon}(z)$ exists for almost all $z\in\Phi$ with respect to the
liouville measure. The function $\vec{\upsilon}:\Phi\to\R^q$ is measurable and
invariant under $\vap:\Phi\to\Phi$. It satisfies
$$
\vec{0}\prec \vec{\upsilon}(\cdot) \prec (|\bo_1\Om|,\dots,|\bo_q\Om|).
$$
We have
\begin{equation}    \label{vec_foot_gain_eq}
\vec{0}\prec \frac{\int_{\Phi}\vec{\upsilon}(z)d\nu(z)}{\nu(\Phi)} \prec(|\bo_1\Om|,\dots,|\bo_q\Om|).
\end{equation}
\end{prop}

\vspace{5mm}

Normalizing $\vec{\upsilon}(\cdot)$, we define the {\em rotation vector}
\begin{equation}    \label{rot_vec_eq}
\vec{\rho}(z)=\left(\frac{\ups_1(z)}{|\bo_1\Om|},\cdots,\frac{\ups_q(z)}{|\bo_q\Om|}\right).
\end{equation}
We reformulate Proposition~~\ref{vec_foot_gain_prop} in terms of the rotation vector.
\begin{prop}    \label{rot_vec_prop}
The rotation vector $\vec{\rho}(z)$ exists for almost all $z\in\Phi$ with respect to
the liouville measure. The vector function $\vec{\rho}:\Phi\to\R^q$ is measurable and
$\vap$-invariant. We have the bounds
$$
\vec{0} \prec \vec{\rho}(\cdot) \prec (1,\dots,1).
$$
The average  rotation vector satisfies
\begin{equation}    \label{av_rot_vec_eq}
\vec{0}\prec \frac{\int_{\Phi}\vec{\rho}(z)d\nu(z)}{\nu(\Phi)} \prec(1,\dots,1).
\end{equation}
\end{prop}

\medskip


%
\begin{thm}   \label{rot_vec_thm}
Let $\Om\subset\RR$ be a bounded, connected, multiply connected domain with a piecewise smooth
boundary. Let $\vap:\Phi\to\Phi$ be the billiard map. Suppose that $\Om$ is $q$-connected, and
let $\vec{\rho}:\Phi\to\R^q$ be the rotation vector function.

\medskip

\noindent The mean rotation vector satisfies
\begin{equation}    \label{mean_rot_vec_eq}
\frac{\int_{\Phi}\vec{\rho}(z)d\nu(z)}{\nu(\Phi)}=\frac12\left(\frac{|\bo_1\Om|}{|\bo\Om|},\cdots,
\frac{|\bo_q\Om|}{|\bo\Om|}\right).
\end{equation}

Suppose that the billiard in $\Phi$ is ergodic. Then
for almost all $z\in\Phi$ we have
$$
\vec{\rho}(z)=\frac12\left(\frac{|\bo_1\Om|}{|\bo\Om|},\cdots,
\frac{|\bo_q\Om|}{|\bo\Om|}\right).
$$
\begin{proof}

Let $\si:\Phi\to\Phi$ and $\tau:\Phi\to\Phi$ be the canonical
involutions. (See the proof of Theorem~\ref{rot_num_thm}.) Both
conjugate $\vap$ with $\vap^{-1}$. It is useful to think of
billiard orbits geometrically as oriented curves in $\Phi$. These
curves are piecewise linear; an orbit $\ga$ is a sequence of
straight segments with endpoints in $\bo\Phi$. In order to
distinguish between orbits of the billiard map $\vap:\Phi\to\Phi$
and the geometric orbits, in what follows we will call the latter
the {\em billiard curves}. The involutions $\si$ and $\tau$ send
billiard curves into themselves, reversing the orientation. The
billiard curves $\si(\ga)$ and $\tau(\ga)$ differ only by a shift
in segment labelling. Since a  finite shift in labelling of
billiard curves does not change the outcome of computing footpoint
gains, we do not distinguish between $\si(\ga)$ and $\tau(\ga)$.
We denote this  infinite billiard curve by $\tga$.

Let $z\in\Phi$ and $\tz\in\Phi$ be a pair of phase points related by a canonical involution.
Let $\ga$ and $\tga$ be the corresponding billiard curves. To compute the numbers $\ups_{\al}(z),\,\ups_{\al}(\tz)$,
we perform the following operations: i) We consider finite subcurves, say $\ga_N$, of $\ga$
(resp. $\tga_N$ of $\tga$) consisting of $N$ segments; ii) We add up the footpoint
gains along $\bo_{\al}\Om$ corresponding to consecutive visits of  $\ga_N$ (resp. $\tga_N$) to the component $\bo_{\al}\Om$;
iii) We divide the result by $N$, and take the limit $N\to\infty$.

Let $\de$ and $\tde$
be the footpoint gains for $\ga$ and $\tga$ corresponding to the same pair of consecutive visits.
Then $\de+\tde=|\bo_{\al}\Om|$. Therefore, the sum $\ups_{\al}(z)+\ups_{\al}(\tz)$ is equal to
the product of $|\bo_{\al}\Om|$ and the
frequency of visits to $\bo_{\al}\Om$ for the infinite billiard curve $\ga$. Integrating this
equality over $\Phi$, using the basic properties of $\nu$ and the
invariance of the liouville measure with respect to $z\mapsto\tz$,
we obtain the identity
\begin{equation}      \label{comp_foot_eq}
\int_{\Phi}\ups_{\al}(z)d\nu(z)=|\bo_{\al}\Om|^2.
\end{equation}

Putting together equations~~\eqref{comp_foot_eq} for $1\le\al\le q$, we arrive at
the vector identity
\begin{equation}      \label{aver_foot_eq}
\frac{\int_{\Phi}\vec{\ups}(z)d\nu(z)}{\nu(\Phi)}=\frac12\left(\frac{|\bo_1\Om|^2}{|\bo\Om|},\cdots,
\frac{|\bo_q\Om|^2}{|\bo\Om|}\right).
\end{equation}

Passing from the footpoint gain $\vec{\ups}(z)$ to the rotation vector, we obtain the identity
equation~~\eqref{mean_rot_vec_eq}. The birkhoff ergodic theorem implies the remaining claim.
\end{proof}
\end{thm}

\medskip

Now we define the concept of rotation number for arbitrary connected domains.

\begin{defin}    \label{rot_num_def}
Let $\Om\subset\RR$ be a bounded, connected, multiply connected, in general, domain  with a piecewise smooth
boundary. Suppose that $\Om$ is $q$-connected, and
let $\vec{\rho}=(\rho_1,\dots,\rho_q):\Phi\to\R^q$ be the rotation vector function.

\medskip

We set $\rho(z)=\rho_1(z)+\cdots+\rho_q(z)$. We call $\rho:\Phi\to\R$ the
{\em rotation number function} for the  domain $\Om$.
\end{defin}

Note that if $\Om$ is simply connected, then $\rho:\Phi\to\R$ in Definition~~\ref{rot_num_def}
coincides with the rotation number that we have defined earlier.
The following is an immediate consequence of Theorem~~\ref{rot_vec_thm} and
Theorem~~\ref{rot_num_thm}.
\begin{corol}   \label{rot_num_cor}
Let $\Om\subset\RR$ be a bounded, connected domain with a piecewise smooth
boundary $\bo\Om$ consisting of an arbitrary number of connected components.
Let $\rho:\Phi\to\R_+$ be the rotation number function for $\Om$.

\noindent Then $\rho$ is a measurable function with values in $[0,1]$. Let
$X\subset[0,1]$ be an arbitrary measurable set. The rotation number function
satisfies the symmetry
$$
\nu(\rho^{-1}(1-X))=\nu(\rho^{-1}(X)).
$$
The mean value of $\rho$ on $\Phi$ with respect
to the liouville measure is $1/2$. If the billiard in $\Om$ is ergodic, then $\rho(z)=1/2$
for $\nu$-almost all $z\in\Phi$.
\end{corol}
\medskip
\section{Examples, amplifications, and open questions}   \label{conclu}
We will now discuss the preceding material in several special but representative cases.
\subsection{Nonergodic billiards in simply connected domains, and related questions}   \label{twist_sub}
Let $\Om\subset\RR$ be a strictly convex domain with a $C^2$ boundary.
Then $\vap:\Phi\to\Phi$ is a $C^1$ twist map. The rotation number $\rho(z)$ introduced in
section~~\ref{simply} coincides with the classical  rotation number for area preserving twist maps
\cite{Ba88}. We point out that equation~~\eqref{mean_rot_eq} does not
have a counterpart for general twist maps. It reflects special features
of the billiard dynamics.

It is an open question whether the billiard in a domain $\Om$ of this kind can be ergodic \cite{Gut03}.
If $\bo\Om$ is of class $C^7$ and its curvature is strictly positive, then the billiard in
$\Om$ is not ergodic. This follows from the results of Lazutkin \cite{La73} and R. Douady \cite{Dou} about
billiard caustics. If the curvature of $\bo\Om$ is not strictly positive, then, by a theorem of Mather \cite{Ma82},
$\Om$ has no caustics. However, this does not imply the ergodicity of the billiard in $\Om$. In particular,
it is conceivable that the essential range of the rotation function is $[0,1]$.

A caustic in $\Om$ yields a topologically nontrivial $\vap$-invariant curve $\Ga\subset\Phi$.
The restriction $\vap|_{\Ga}:\Ga\to\Ga$ is a circle homeomorphism. Let $\rho(\Ga)$ be the Poincare
rotation number for $\vap|_{\Ga}$. Then $\rho(z)=\rho(\Ga)$ for all $z\in\Ga$. This raises the question:
for which $\Om$ there is an open region $X\subset\Phi$ foliated by invariant curves? By a conjecture
attributed to G.D. Birkhoff \cite{Gut03}, the only regions satisfying this property are ellipses.
The question is open except for an important special case when $X=\Phi$. Then, by a theorem of Bialy \cite{Bi93},
$\Om$ is a disc. See also \cite{Woj94} for the relevant material.

Let $\Ga\subset\Phi$ be an invariant curve. What can we say about the rotation number $\rho(\Ga)$?
Kolodziej \cite{Ko85} explicitly computed $\rho(\Ga)$ in the case when $\Om$ is an ellipse.
Although the answer is in terms of elliptic functions, Kolodziej's approach is geometric.
Namely, the work \cite{Ko85} crucially uses the invariant measure for a natural geometric transformation associated with
a pair of nested circles.

\medskip
\subsection{Ergodic billiards in simply connected domains}   \label{erg_sub}
By Theorem~~\ref{rot_num_thm}, ergodicity of the billiard in $\Om$
insures that $\rho=1/2$ almost everywhere. This raises the
question: For which simply connected $\Om$ the billiard is
ergodic? There is a vast literature on the subject; see
\cite{LW95} and the references there. As a rule, ergodic domains
are not convex and their boundaries have singular points. Although
there are examples of convex ergodic $\Om$ (e. g., the {\em
stadium}) but $\bo\Om$ in these examples has low regularity. In
particular, the stadium is only $C^1$ and not strictly convex.

Ergodicity of the billiard domains discussed in  \cite{LW95} is
associated with the {\em hyperbolicity}. There are broad geometric
conditions on $\bo\Om$ that ensure hyperbolicity \cite{Woj86}.
However, a hyperbolic billiard domain may be non-ergodic; there
are explicit examples of hyperbolic $\Om$ with several ergodic
components \cite{Woj86}. It would be interesting to compute the
rotation numbers for  ergodic components in these examples.

There is a class of (possibly) ergodic simply connected domains
with parabolic billiard dynamics. Let $P\subset\RR$ be a simply
connected polygon. It is {\em rational} if all of its angles are
commensurable with $\pi$. Otherwise $P$ is {\em irrational}. If
$P$ is rational, then the direction of a tangent vector in $P$ is
preserved modulo the action of a finite group \cite{Gut84}. In
view of this, the billiard phase space is foliated by a
one-parameter family of invariant subsets, hence the billiard in a
rational polygon $P$ is not ergodic. If $P$ is irrational, then
there is no obvious obstruction to ergodicity. Presently there are
no examples of nonergodic irrational polygons \cite{Gut03}. A
theorem of Vorobets \cite{Vo97} says that if the angles of an
irrational polygon $P$ are extremely fast approximated by
$\pi$-rational numbers, then $P$ is ergodic. This result yields
explicit examples of simply connected (even convex) ergodic
polygons. By Theorem~~\ref{rot_num_thm}, for such $P$ we have
$\rho=1/2$ almost everywhere. It would be interesting to compute
the function $\rho(z)$ for rational polygons. The case of {\em
arithmetic} polygons \cite{GJ00} seems especially intriguing.

\medskip
\subsection{Billiards in multiply connected domains}   \label{multi_sub}
Multiply connected billiard domains arise naturally in several
contexts. We will discuss one of them. Let $\ga\subset\RR$ be a
closed, convex curve. The {\em string construction} \cite{GK95}
yields a one-parameter family $X=X(\ga,l)$ of convex domains
containing $\ga$ in their interior. The parameter $l>|\ga|$ is the
string length. The family $X(\ga,l)$ consists of planar domains
having $\ga$ as a billiard caustic \cite{GK95}. Let $\Phi(l)$ be
the phase space for the billiard map in $X(\ga,l)$; let
$\Ga(l)\subset\Phi(l)$ be the invariant curve corresponding to
$\ga$; let $\rho(l)$ be its rotation number. The monotonic
function $l\mapsto\rho(l)$ yields the so called {\em phase
locking} phenomenon \cite{GKn96}. The case when $\ga$ is a polygon
is especially interesting. Then the domains  $X=X(\ga,l)$ exhibit
the phenomenon of {\em instability of the boundary} \cite{Hub87},
\cite{GK95}. This means that near $\bo X$ there is an open region
in $X$ which is {\em free of caustics}. It is plausible that
$\ga\subset X(\ga,l)$ is the {\em last caustic} \cite{GKn96}.

Let $\Om=X(\ga,l)\setminus\inter(\ga)$. Thus $\Om=\Om(\ga,l)$ is a
multiply connected domain. It would be instructive to investigate
the rotation vector for the billiard in $\Om$. See
section~~\ref{multi}. The case when $\ga$ is a circle is
especially simple.  Then the regions $\Om(\ga,l)$ are the
rotationally symmetric annuli. The billiard in a rotationally
symmetric annulus is completely integrable; concentric circles
yield invariant curves that foliate the phase space. When $\ga$ is
an ellipse, the situation is more complicated, but the billiard in
question is still completely integrable. There are two kinds of
invariant curves: confocal ellipses and confocal hyperbolas. It is
plausible that these are the only multiply connected domains with
completely integrable billiard dynamics. This is a natural
extension of the conjecture attributed to G.D. Birkhoff. See
section~~\ref{twist_sub}.

The above multiply connected domains can be described as follows:
$\Om$ is the annulus between a closed convex curve $\ga$ and the
boundary of a billiard table, say $X$, having $\ga$ as a billiard
caustic. An interesting class of multiply connected domains arises
when we perturb this situation in a particular way. Namely, $\Om$
is obtained by moving $\ga$ inside $X$ by a small planar isometry.
In the simplest case, when $\ga$ is a circle, we obtain annular
regions between two non-concentric circles. We will refer to these
regions as {\em asymmetric annuli}. Let $\Om$ be an asymmetric
annulus; let $\Phi$ be the corresponding billiard phase space.
There is an obvious invariant region $\Phi_0\subset\Phi$ foliated
by invariant curves. The restriction of billiard dynamics to
$\Phi_0$ is completely integrable. The analysis of hyperbolicity
and ergodicity in the complement $\Phi\setminus\Phi_0$ by
currently available methods is inconclusive. However, numerical
simulations indicate that the billiard dynamics in
$\Phi\setminus\Phi_0$ is chaotic.

\end{document}